\documentclass{amsart}

\usepackage{amsmath,amssymb}
\usepackage[initials]{amsrefs}

\newtheorem{theorem}{Theorem}[section]

\newtheorem{lemma}[theorem]{Lemma}
\newtheorem{proposition}[theorem]{Proposition}

\theoremstyle{definition}

\theoremstyle{remark}
\newtheorem*{remark}{Remark}

\newcommand{\F}{\mathbb{F}}

\newcommand{\ord}{\mathrm{ord}}

\newcommand{\PI}{\mathcal{P}}
\newcommand{\Li}{\mathcal{L}}
\newcommand{\T}{\mathcal{T}}
\newcommand{\X}{\mathcal{X}}

\newcommand{\Z}{\mathbb{Z}}
\newcommand{\N}{\mathcal{N}}

\title[Extensions with the line property]{Finite field extensions with the line or translate property for $\lowercase{r}$-primitive elements}

\author{Stephen D. Cohen}\thanks{The first author is Emeritus Professor of Number Theory, University of Glasgow}
\address{6 Bracken Road, Portlethen, Aberdeen AB12 4TA, Scotland, UK}
\email{Stephen.Cohen@glasgow.ac.uk}

\author{Giorgos Kapetanakis}
\address{Department of Mathematics and Applied Mathematics, University of Crete, Voutes Campus, 70013 Heraklion, Greece}
\email{gnkapet@gmail.com}

\subjclass[2010]{11T30 (Primary), 11T06 (Secondary)}

\keywords{Primitive element, high order element, line property, translate property}

\begin{document}

\date{\today}

\begin{abstract}
Let $r,n>1$ be integers and $q$ be any  prime power $q$ such that $r\mid q^n-1$.  We say that the extension  $\F_{q^n}/\F_q$ possesses the line property for $r$-primitive elements  property if,  for every $\alpha,\theta\in\F_{q^n}^*$, such that $\F_{q^n}=\F_q(\theta)$, there exists some $x\in\F_q$, such that $\alpha(\theta+x)$ has multiplicative order $(q^n-1)/r$. We prove that, for sufficiently large prime powers $q$, $\F_{q^n}/\F_q$ possesses the line property for $r$-primitive elements. We also discuss the (weaker) translate property for extensions.
\end{abstract}

\maketitle

\section{Introduction}
Let $q$ be a prime power and $n\geq 2$ an integer. We denote by $\F_q$ the finite field of $q$ elements and by $\F_{q^n}$ its extension of degree $n$.
It is well-known that the multiplicative group $\F_{q^n}^*$ is cyclic; its  generators are called  \emph{primitive elements}. The theoretical importance of primitive elements is complemented by their numerous applications in practical areas such as cryptography.

In addition to their theoretical interest, elements of $\F_{q^n}^*$ that have high order, without necessarily being primitive, are of great practical interest because in several applications they may replace primitive elements.   Accordingly, recently researchers have worked on the the effective construction of such high order elements, \cites{gao99,martinezreis16,popovych13}, since that  of  primitive elements themselves remains an open problem.

With that in mind, we call an element of order $(q^n-1)/r$, where $r\mid q^n-1$, \emph{$r$-primitive}, i.e., the primitive elements are exactly the $1$-primitive elements. In this line of work, the existence of $2$-primitive elements that also possess other desirable properties has been recently considered \cites{cohenkapetanakis19,kapetanakisreis18}.

We call some $\theta\in\F_{q^n}$ a \emph{generator} of the extension $\F_{q^n}/\F_q$ if $\F_{q^n} = \F_q(\theta)$ and, if $\theta$ is a generator of $\F_{q^n}/\F_q$, we call the set \[ \T_\theta := \{ \theta +x \, :\, x\in\F_q\} \] the \emph{set of translates} of $\theta$ over $\F_q$ and every element of this set a \emph{translate} of $\theta$ over $\F_q$.
We say that an extension $\F_{q^n} / \F_q$ possesses the \emph{translate property for $r$-primitive elements}, if every set of translates contains an $r$-primitive element. In particular, for $r=1$ we simply call it the \emph{translate property}.
A classical result in the study of primitive elements is the following.
\begin{theorem}[Carlitz-Davenport]\label{thm:ca-da}
Let $n$ be an integer. There exist some $T_1(n)$ such that for every prime power $q> T_1(n)$, the extension $\F_{q^n}/\F_q$ possesses the translate property.
\end{theorem}
The above was first proved by Davenport \cite{davenport37}, for prime $q$, while Carlitz \cite{carlitz53a} extended it to the stated form. Interest in this problem was renewed by recent applications of the translate property in semifield primitivity, \cites{kapetanakislavrauw18,rua15,rua17}.

Let $\theta$ be a generator of the extension $\F_{q^n}/\F_q$ and take some $\alpha\in\F_{q^n}^*$. We call the set \[ \Li_{\alpha,\theta} := \{ \alpha(\theta+x) \, : \, x\in\F_q \} \] the \emph{line} of $\alpha$ and $\theta$ over $\F_q$. An extension $\F_{q^n}/\F_q$ is said to possess the \emph{line property for $r$-primitive elements} if every line of this extension contains an $r$-primitive element. When $r=1$, we refer to this property as the \emph{line property}.
A natural generalization of Theorem~\ref{thm:ca-da} is the following, \cite{cohen10}*{Corollary~2.4}.
\begin{theorem}[Cohen] \label{thm:line1}
Let $n$ be an integer. There exist some $L_1(n)$ such that for every prime power $q> L_1(n)$, the extension $\F_{q^n}/\F_q$ possesses the line property.
\end{theorem}
It is clear that all the sets of translates are actually lines (where $\alpha=1$), i.e., the line property implies the translate property. Thus, $L_1(n) \geq T_1(n)$.
In this work, we extend Theorems~\ref{thm:ca-da} and \ref{thm:line1} to $r$-primitive elements, by proving the following.
\begin{theorem}\label{thm:our_ca-da}
Let $n$ and $r$ be integers. There exist some $L_r(n)$ such that for every prime power $q> L_r(n)$, with the property $r\mid q^n-1$, the extension $\F_{q^n}/\F_q$ possesses the line property for $r$-primitive elements. If we confine ourselves to the translate property for $r$-primitive elements, the same is true for some $T_r(n)\leq L_r(n)$.
\end{theorem}
A natural, but apparently challenging, related question is identifying the exact value of the numbers $T_1(n)$ and $L_1(n)$ for given $n$.  Indeed, only  a handful of these are known.  In particular, the first author, in \cite{cohen83}, proved that $T_1(2) = L_1(2) = 1$ and, in \cite{cohen09}, that $T_1(3) = 37$. Bailey et al.~\cite{baileycohensutherlandtrudgian19} proved that $L_1(3) = 37$ and estimated $T_1(4) \leq L_1(4) \leq 102829$. By means of more specialized theoretical and computational techniques, the authors have now proved that $T_2(2)=L_2(2)=41$.   The details will be given in a further article.
\section{Preliminaries}\label{sec:preliminaries}
We begin by introducing the notion of freeness. Let $m\mid q^n-1$, an element $\xi \in \F_{q^n}^*$ is \emph{$m$-free} if $\xi = \zeta^d$ for some $d\mid m$ and $\zeta\in\F_{q^n}^*$ implies $d=1$. It is clear that primitive elements are exactly those that are $q_0$-free, where $q_0$ is the square-free part of $q^n-1$. It is also evident that there is some relation between $m$-freeness and multiplicative order.
\begin{lemma}[\cite{huczynskamullenpanariothomson13}, Proposition~5.3]\label{lemma:m-free}
If $m\mid q^n-1$ then $\xi\in\F_{q^n}^*$ is $m$-free if and only if $\gcd\left( m,\frac{q^n-1}{\ord\xi} \right)=1$.
\end{lemma}

Throughout this work, a \emph{character} is a multiplicative character of $\F_{q^n}^*$, while we denote by $\chi_0$ the trivial multiplicative character.
Vinogradov's formula yields an expression for the characteristic function of $m$-free elements in terms of multiplicative characters, namely:
\begin{equation}\label{eq:omega}
\Omega_m(x) := \theta(m) \sum_{d\mid m} \frac{\mu(d)}{\phi(d)} \sum_{\ord\chi = d} \chi(x) ,
\end{equation}
where $\mu$ stands for the M\"{o}bius function, $\phi$ for the Euler function, $\theta(m) := \phi(m)/m$ and the inner sum suns through multiplicative characters of order $d$. Furthermore, a direct consequence of the orthogonality relations is that the characteristic function for the elements of $\F_{q^n}^*$ that are $k$-th powers, where $k\mid q^n-1$, can be written as
\begin{equation}\label{eq:w}
w_k (x) := \frac{1}{k} \sum_{d\mid k} \sum_{\ord\chi=d} \chi(x) .
\end{equation}

We will use character sums to establish our results. The following estimate is a direct consequence of the main result of \cite{katz89} and, notably, it is one of the few non-trivial character sum estimates not relying on Weil's results \cite{weil48b}.
\begin{proposition}[Katz]\label{prop:katz}
Let $\theta\in\F_{q^n}$ be such that $\F_{q^n} = \F_q(\theta)$ and $\chi$ a non-trivial character. Then
\[
\left| \sum_{x\in\F_q} \chi (\theta+x) \right| \leq (n-1) \sqrt{q} .
\]
\end{proposition}
Let $d(R)$ be the number of divisors of $R$. The result below provides an asymptotic estimate for this function.
\begin{proposition}[\cite{apostol76}, p.~296] \label{prop:apostol}
For every $\delta>0$, $d(n) = o(n^\delta)$, where $o$ signifies the little-o notation.
\end{proposition}
\section{Characterization of $r$-primitive elements}\label{sec:char}
From now on, fix the positive integers $r$ and $n$ and let $q$ be some prime power such that $r\mid q^n-1$. Let $\Gamma$ be the characteristic function for $r$-primitive elements of $\F_{q^n}$, that is,
\[
\Gamma (x) := \begin{cases} 1, & x\text{ is $r$-primitve}, \\ 0, & \text{otherwise.} \end{cases},\ x\in\F_{q^n}^*.
\]
The aim of this section is to express $\Gamma(x)$ in a convenient way, using characters.

Let $\PI$ be the set of distinct primes dividing $q^n-1$. It follows that $q^n-1 = \prod_{p\in\PI} p^{a_p}$, where $a_p\geq 1$ for all $p\in\PI$. Additionally, we have that $r=\prod_{p\in\PI} p^{b_p}$, where for every $p\in\PI$,   $0\leq b_p\leq a_p$.

We partition $\PI$ as follows:
\begin{align*}
\PI_s & := \{ p\in\PI\, : \, a_p=b_p > 0 \} ,\\
\PI_t & := \{ p\in\PI\, : \, a_p>b_p > 0 \} ,\\
\PI_u & := \{ p\in\PI\, : \, a_p>b_p = 0 \} .
\end{align*}
It is clear that the above sets are pairwise disjoint and that $\PI_s \cup \PI_t\cup\PI_u=\PI$.
Further, set
\[ s:=\prod_{p\in\PI_s} p^{b_p},\, t:=\prod_{p\in\PI_t} p^{b_p}\text{ and } u:= \prod_{p\in\PI_u} p . \]
It is straightforward to check that $r=st$ and that $u$ is the radical of the   part of $q^n-1$ that is relatively prime with $r$.

Lemma~\ref{lemma:m-free}, implies that the set of $u$-free elements, contains all the $\sigma$-primitive elements, where
\begin{equation} \label{eq:sigma}
\sigma = \prod_{p\in\PI_s\cup\PI_t} p^{\sigma_p} ,
\end{equation}
for some $0\leq \sigma_p\leq a_p$. In addition, the $u$-free elements that are $r$-th powers are the $\sigma$-primitive elements with $\sigma$ as in \eqref{eq:sigma} with $b_p \leq \sigma_p \leq a_p$. 
Next, let $\PI_t = \{ p_1 ,\ldots p_k \}$ and for $i=1,\ldots , k$ set $e_i := p_i^{b_{p_i}}$ and $f_i := p_i^{b_{p_i}+1}$. Note that $b_{p_i}+1\leq a_{p_i}$. Now from the set of $u$-free elements, that are also $r$-th powers, exclude those that are not $f_i$-th powers for every $i=1,\ldots,k$. We are left with exactly the $\sigma$-primitive elements, where $\sigma$ is as in \eqref{eq:sigma}, with
\[
\begin{cases}
b_p \leq \sigma_p \leq a_p = b_p , & \text{if } p\in\PI_s \\
b_p \leq \sigma_p < b_p+1 \leq a_p , & \text{if } p\in\PI_t .
\end{cases}
\]
In particular, in any case $\sigma_p = b_p$, that is, $\sigma = r$.

In other words, with the notation of Section~\ref{sec:preliminaries}, the characteristic function for $r$-primitive elements of $\F_{q^n}^*$ can be expressed as
\begin{align} 
\Gamma(x) & = \Omega_{u}(x) w_{r}(x) \prod_{i=1}^k (1-w_{f_i}(x)) \nonumber \\
 & = \Omega_{u}(x) w_{s}(x) \prod_{i=1}^k w_{e_i}(x)(1-w_{f_i}(x)) , \label{eq:Omega1}
\end{align}
where $x\in\F_{q^n}^*$. Moreover, for every $i=1,\ldots,k$, notice that since $e_i\mid f_i$, we have that an $f_i$-th power is also an $e_i$-th power, i.e., for $x\in\F_{q^n}^*$, $w_{e_i}(x)w_{f_i}(x) = w_{f_i}(x)$. Thus \eqref{eq:Omega1} yields
\begin{equation} \label{eq:Omega2}
\Gamma(x) = \Omega_{u}(x) w_{s}(x) \prod_{i=1}^k (w_{e_i}(x)-w_{f_i}(x)).
\end{equation}
Next, recall that, for $i=1,\ldots,k$, $e_i=p_i^{b_{p_i}}$ and $f_i=p_i^{b_{p_i}+1} = e_ip_i$. It follows that for every $x\in\F_{q^n}^*$
\begin{align*}
w_{e_i}(x)-w_{f_i}(x) & = \frac{1}{e_i} \sum_{d\mid e_i} \sum_{\ord\chi=d} \chi(x) - \frac{1}{f_i} \sum_{d\mid f_i} \sum_{\ord\chi=d} \chi(x) \\
 & = \frac{1}{e_i} \sum_{d\mid f_i} \sum_{\ord\chi = d} \ell_{i,d} \chi(x) ,
\end{align*}
where, for $d\mid f_i$,
\[ 
\ell_{i,d} := \begin{cases} 1 - 1/p_i , & \text{if } d\neq f_i , \\ -1/p_i , & \text{if } d=f_i . \end{cases}
\]
Finally, we insert the above and the expressions \eqref{eq:omega} and \eqref{eq:w} into \eqref{eq:Omega2}, and  obtain
\begin{equation} \label{eq:Omega_expanded}
\Gamma(x) = \frac{\theta(u)}{r} \sum_{d_1\mid u} \sum_{d_2\mid s} \sum_{\delta_1\mid f_1} \cdots \sum_{\delta_\lambda\mid f_k} \frac{\mu(d_1)}{\phi(d_1)} \ell_{1,\delta_1} \cdots \ell_{k,\delta_k} \sum_{\substack{\ord\chi_j=d_j \\ \ord\psi_i=\delta_i}} (\chi_1\chi_2\psi_1\cdots\psi_k)(x) ,
\end{equation}
where $x\in\F_{q^n}^*$ and $(\chi_1\chi_2\psi_1\cdots\psi_\lambda)$ stands for the product of the corresponding characters, itself a character.
\section{Proof of Theorem~\ref{thm:our_ca-da}}\label{sec:main}
Fix some $\alpha,\theta\in\F_{q^n}^*$, such that $\F_{q^n}=\F_q(\theta)$. Let $\N(\theta,\alpha)$ be the number of $r$-primitive elements of the form $\alpha(\theta+x)$, where $x\in\F_q$. It suffices to show that
\[
\N(\theta,\alpha) = \sum_{x\in\F_q} \Gamma(\alpha(\theta+x)) \neq 0 .
\]
With \eqref{eq:Omega_expanded} in mind, we have that
\begin{equation}\label{eq:N1}
\frac{\N(\theta,\alpha)}{\theta(u)} = \frac{1}{r} \sum_{\substack{ d_1\mid u,\, d_2\mid s, \\ \delta_1\mid f_1 ,\ldots , \delta_k\mid f_k}} \frac{\mu(d_1)}{\phi(d_1)} \ell_{1,\delta_1} \cdots \ell_{k,\delta_k} \sum_{\substack{\ord\chi_j=d_j \\ \ord\psi_i=\delta_i}} \X_{\alpha,\theta}(\chi_1,\chi_2,\psi_1,\ldots ,\psi_k) ,
\end{equation}
where
\[
\X_{\alpha,\theta}(\chi_1,\chi_2,\psi_1,\ldots ,\psi_k) := \sum_{x\in\F_q} (\chi_1\chi_2\psi_1\cdots\psi_k)(\alpha(\theta+x)) .
\]
In addition, notice that the orders of all the factors of the character product $(\chi_1\chi_2\psi_1\cdots\psi_k)$ are relatively prime. Hence the product itself is trivial if and only if all its factors are trivial. With this in mind, Proposition~\ref{prop:katz} implies that, unless all the characters $\chi_1$, $\chi_2$, $\psi_1$,\ldots,$\psi_k$ are trivial,
\[
| \X_{\alpha,\theta}(\chi_1,\chi_2,\psi_1,\ldots ,\psi_k) | \leq \sqrt{q} ,
\]
while it is clear that
\[
\X_{\alpha,\theta}(\chi_0,\chi_0,\chi_0,\ldots ,\chi_0) = q.
\]
In \eqref{eq:N1}, we separate the term that corresponds to $d_1=d_2=\delta_1=\ldots=\delta_k = 1$ and, with the above in mind, we obtain
\begin{equation}\label{eq:N2}
\left| \frac{\N(\theta,\alpha)}{\theta(u)} - \frac{q}{r} \cdot \ell_{1,1} \cdots \ell_{k,1} \right| \leq 
\frac{1}{r} \sum_{\substack{ d_1\mid u,\, d_2\mid s, \\ \delta_1\mid f_1 ,\ldots , \delta_k\mid f_k \\ \text{not all equal to }1}} \frac{|\ell_{1,\delta_1} \cdots \ell_{k,\delta_k}|}{\phi(d_1)} \sum_{\substack{\ord\chi_j=d_j \\ \ord\psi_i=\delta_i}} \sqrt{q} .
\end{equation}
Notice that for all $1\leq i\leq k$, $|\ell_{i,\delta_i}|\leq \ell_{i,1}$. It follows from \eqref{eq:N2} that $\N(\theta,\alpha)\neq 0$ if
\[
q > \sum_{\substack{ d_1\mid u,\, d_2\mid s, \\ \delta_1\mid f_1 ,\ldots , \delta_k\mid f_k}} \frac{1}{\phi(d_1)} \sum_{\substack{\ord\chi_j=d_j \\ \ord\psi_i=\delta_i}} \sqrt{q} .
\]
Furthermore, it is well-known that, for every $d\mid q^n-1$, there exist exactly $\phi(d)$ characters of order $d$.  Hence the latter condition can be also written as
\begin{equation} \label{eq:cond_main}
q >   s f_1\cdots f_k \cdot d(u) \cdot \sqrt{q},
\end{equation}
where we recall that $d(m)$ stands for the number of divisors of $m\in\Z$. Now, observe that $u \mid q^{n}-1$, wherefore Proposition~\ref{prop:apostol} implies that
\[ d(u) \leq d(q^{n}-1)= o(q^{1/4}). \] Further, observe that
\[ sf_1\cdots f_k \leq A_r := \prod_{p\in\PI_s\cup \PI_t} p_i^{b_i+1} , \]
where the left side of the above inequality depends solely on $r$. It follows that, for $q$ large enough, \eqref{eq:cond_main} holds.  Hence $\N(\theta,\alpha)\neq 0$.

The proof of the first statement of Theorem~\ref{thm:our_ca-da} is now complete, while the second statement (about $T_r(n)$) follows immediately from the fact that all the sets of translates of an extension are simultaneously lines of this extension.
\begin{remark}
The reader will note that, in the statement of Theorem~\ref{thm:our_ca-da}, $n$ and $r$ are fixed integers and $q$ is any (sufficiently large) prime power such that  $r\mid q^n-1$.  We are indebted to the referee for the observation that it would be possible to allow the integer $r$  with this property to vary as a (suitably small) function of $q$. To avoid unnecessary technical complication we have refrained from a precise formulation of the ensuing result at this  juncture.
\end{remark}
\section*{Acknowledgments}
We are grateful to Michel Lavrauw  whose suggested terminology we have adopted.
%
%
%
%
%

%
\end{document}